\documentclass{article}

\usepackage{amssymb}
\usepackage{amsmath}

\usepackage[dvips]{graphicx}

\begin{document}


\begin{center}
\noindent \textbf{On the paper: Numerical radius preserving linear maps on Banach algebras}                                                                                                                                                                                     
\end{center}

\noindent \textbf{}

\begin{center}
M. El Azhari
\end{center}

\noindent \textbf{ } 

\noindent \textbf{Abstract.} We give an example of a unital commutative complex Banach algebra having a normalized state which is not a spectral state and admitting an extreme normalized state which is not multiplicative. This disproves two results by Golfarshchi and Khalilzadeh.

\noindent \textbf{}

\noindent \textit{Keywords and phrases:} Banach algebra, regular norm, normalized state, spectral state.

\noindent \textbf{}

\noindent \textit{(2010) Mathematics Subject Classification:} 46H05.

\noindent \textbf{}

\noindent \textbf{} 

\noindent \textbf{1. Preliminaries}

\noindent \textbf{} 

\noindent \textbf{} Let $(A,\Vert .\Vert)$ be a complex normed algebra with an identity $e$ such that $\Vert e\Vert = 1.$ Let $D(A,e)=\lbrace f\in A^{'}: f(e)=\Vert f\Vert = 1\rbrace,$ where $A^{'}$ is the dual space of $A.$ The elements of $D(A,e)$ are called normalized states on $A.$ For $a\in A,$ let $V(A,a)=\lbrace f(a): f\in D(A,e)\rbrace,\: V(A,a)$ is called the numerical range of $a.$ Let $sp(a) $ be the spectrum of $a\in A,$ and let $co(sp(a))$ be the convex hull of $sp(a).$  We say that a linear functional $ f $ on $ A $ is a spectral state if $f(a)\in co(sp(a))$ for all $a\in A.$ We denote by $M(A)$ the set of all non-zero continuous multiplicative linear functionals on $A.$ 

\noindent \textbf{} 

\noindent \textbf{2. Result}

\noindent \textbf{}
 
\noindent \textbf{2.1. Counterexample}

\noindent \textbf{}
 
\noindent \textbf{} Golfarshchi and Khalilzadeh proved the following results [4]: 

\noindent \textbf{}

\noindent \textbf{[4, Theorem 2].} Let $A$ be a unital complex Banach algebra, and let $f$ be a linear functional on $A.$ Then $f$ is a normalized state on $A$ if and only if $f(a)\in co(sp(a))$ for all $a\in A.$

\noindent \textbf{}

\noindent \textbf{[4, Theorem 3].} Let $A$ be a unital commutative complex Banach algebra. Then each  extreme normalized state on $A$ is multiplicative.

\noindent \textbf{}

\noindent \textbf{} Here we give a counterexample disproving the above results. We also remark that Theorems 5 and 6 [4] are called into question since the authors used Theorem 3 [4] to prove these results.

\noindent \textbf{}

\noindent \textbf{} Let $(A,\Vert .\Vert)$ be a non-zero commutative radical complex Banach algebra [6, p.316]. Let $A_{e}=\lbrace a+\lambda e: a\in A,\, \lambda\in C\rbrace$ be the unitization of $A$ with the identity $e,$ and the norm $\Vert a+\lambda e\Vert_{1}=\Vert a\Vert + \vert\lambda\vert$ for all $a+\lambda e\in A_{e}.\: (A_{e},\Vert .\Vert_{1})$ is a unital commutative complex Banach algebra, and $M(A_{e})=\lbrace \varphi_{\infty}\rbrace,$ where $\varphi_{\infty}$ is the continuous multiplicative linear functional on $A_{e}$ defined by $\varphi_{\infty}(a+\lambda e)=\lambda$ for all $a+\lambda e\in A_{e}.$  

\noindent \textbf{(1)} Let $a$ be a non-zero element of $A,\: V(A_{e},a)=\lbrace z\in C: \vert z\vert\leq\Vert a\Vert\rbrace$ by [2, Remark 3.8], and $sp(a)=\lbrace\varphi_{\infty}(a)\rbrace =\lbrace 0\rbrace,$ hence $co(sp(a))=\lbrace 0\rbrace$ is strictly included in  $V(A_{e},a)$ since $\Vert a\Vert\neq 0.$ Therefore the direct implication of [4, Theorem 2] does not hold.

\noindent \textbf{(2)} By [1, lemma 1.10.3], $D(A_{e},e)$ is a non-empty weak* compact convex subset of $A_{e}^{'},$ then $ext (D(A_{e},e))$ is a non-empty set. Assume that each extreme normalized state on $A_{e}$ is multiplicative, then $ext(D(A_{e},e))=\lbrace \varphi_{\infty}\rbrace.$ Let $a$ be a non-zero element of $A,$ by [1, Corollary 1.10.15] there exists $f\in D(A_{e},e)$ such that $f(a)\neq 0=\varphi_{\infty}(a).$ Therefore $\overline{co}(ext(D(A_{e},e)))=\lbrace\varphi_{\infty}\rbrace$ is strictly included in $D(A_{e},e),$ which contradicts the Krein-Milman Theorem. This shows that [4, theorem 3] is not valid.

\noindent \textbf{}

\noindent \textbf{2.2. Regular norm and the operator seminorm} 

\noindent \textbf{} 

\noindent \textbf{} Let $(A,\Vert .\Vert)$ be a non-unital complex Banach algebra, and let  $A_{e}=\lbrace a+\lambda e: a\in A,\, \lambda\in C\rbrace$ be the unitization of $A$ with the identity $e.$ Let $\Vert a+\lambda e\Vert_{op}=\sup \lbrace\Vert (a+\lambda e)x\Vert,\,\Vert x(a+\lambda e)\Vert: x\in A, \Vert x\Vert\leq 1\rbrace$ for all $a+\lambda e\in A_{e},\: \Vert .\Vert_{op}$ is an algebra seminorm on $A_{e}.$ We say that $\Vert .\Vert$ is regular if $\Vert .\Vert_{op}=\Vert .\Vert$ on $A.$ If $\Vert .\Vert$ is regular, it is well known that $(A_{e},\Vert .\Vert_{op})$ is a complex Banach algebra. The following question was asked [3]: If $(A_{e},\Vert .\Vert_{op})$ is a complex Banach algebra, is the norm $\Vert .\Vert$ regular ?

\noindent \textbf{} Orenstein tried to give an answer to this question in the commutative case [5], but his proof is not correct since it is essentially based on the direct implication of [4, Theorem 2].

\noindent \textbf{}

\noindent \textbf{3. Conclusion}

\noindent \textbf{} 

\noindent \textbf{} In this note, we show that Theorems 2 and 3 [4] are false by giving a counterexample. We also remark that Theorems 5 and 6 [4] and Theorem 1.1 [5] are called into question since the authors used Theorems 2 or 3 [4] to prove these results.

\noindent \textbf{} 

\noindent \textbf{} 

\noindent \textbf{References}

\noindent \textbf{}

\noindent \textbf{}[1] F. F. Bonsall and J. Duncan, Complete normed algebras, New York: Springer Verlag 1973.

\noindent \textbf{}[2] A. K. Gaur and T. Husain, Spatial numerical ranges of elements of Banach algebras, International Journal of Mathematics and Mathematical Sciences, 12(4)(1989), 633-640.

\noindent \textbf{}[3] A. K. Gaur and Z. V. Kov\'{a}\v{r}\'{i}k, Norms, states and numerical ranges on direct sums, Analysis, 11(2-3)(1991), 155-164.
 
\noindent \textbf{}[4] F. Golfarshchi and A. A. Khalilzadeh, Numerical radius preserving linear maps on Banach algebras, International Journal of Pure and Applied Mathematics, 88(2)(2013), 233-238.
 
\noindent \textbf{}[5] A. Orenstein, Regular norm and the operator seminorm on a non-unital complex commutative Banach algebra, arXiv:1410.8790v2 [math.FA] 11 Dec 2015.

\noindent \textbf{}[6] C. E. Rickart, General theory of Banach algebras, New York: Van Nostrand 1960.

\noindent \textbf{}

\noindent \textbf{}

\noindent \textbf{}Ecole Normale Sup\'{e}rieure 
  
\noindent \textbf{}Avenue Oued Akreuch

\noindent \textbf{}Takaddoum, BP 5118, Rabat

\noindent \textbf{}Morocco 

\noindent \textbf{} 

\noindent \textbf{}E-mail: mohammed.elazhari@yahoo.fr

\end{document}